
\input amstex
\documentstyle{amsppt}
\pagewidth{6.4in}
\vsize8.5in
\parindent=6mm
\parskip=3pt
\baselineskip=14pt
\tolerance=10000
\hbadness=500
\magnification=1000
\rightheadtext {Spherical maximal operators on
step two nilpotent Lie groups}
\loadbold
\topmatter
\title
Singular spherical maximal operators on a class of step two
nilpotent Lie groups
\endtitle
\author Detlef M\"uller \qquad \qquad \qquad \qquad \qquad
Andreas Seeger \endauthor
\thanks The second author was supported in part by the National Science
Foundation.
\endthanks
\address
Mathematisches Seminar,
Christian-Albrechts-Universit\"at zu Kiel,  Ludewig-Meyn-Str. 4,
24098 Kiel, Germany
\endaddress
\email mueller\@math.uni-kiel.de\endemail
\address
Department of Mathematics,
480 Lincoln Drive, University of Wisconsin, Madison, WI 53706, USA
\endaddress
\email seeger\@math.wisc.edu \endemail
\abstract Let $H^n\cong \Bbb R^{2n}\ltimes \Bbb R$ be the Heisenberg group and let $\mu_t$
be the normalized  surface  measure for the sphere of radius $t$ in
$\Bbb R^{2n}$. Consider
the maximal function defined by $Mf=\sup_{t>0} |f*\mu_t|$.
We prove for $n\ge 2$ that
$M$ defines an operator bounded on $L^p(H^n)$ provided that
$p>2n/(2n-1)$. This improves an earlier result by Nevo and Thangavelu, and the range
for $L^p$ boundedness is optimal.
We also extend the result to  a more general setting of surfaces and to 
groups satisfying a nondegeneracy condition; these include 
the  groups of Heisenberg type.
\endabstract
\subjclass 42B25, 22E25, 43A80
\endsubjclass
\keywords spherical maximal operators, Heisenberg groups,  
step two nilpotent groups,
oscillatory integral operators, fold singularities
\endkeywords\endtopmatter


\def\rk{m}


\def\Ga{\Gamma}

\def\L{{\Cal L}}

\def\tfg{{\widetilde {\frak g}}}

\def\ga{{\gamma}}


\define\rank{{\text{\rm rank }}}

\redefine\ker{{\text{\rm Ker }}}

\define\supp{{\text{\rm supp}}}

\define\inn#1#2{\langle#1,#2\rangle}

\define\lcontr{\rfloor}
\define\lco#1#2{{#1}\lcontr{#2}}
\define\lcoi#1#2{\imath({#1}){#2}}
\define\rco#1#2{{#1}\rcontr{#2}}

\define\bin#1#2{{\pmatrix {#1}\\{#2}\endpmatrix}}
\define\meas{{\text{\rm meas}}}

\define\lc{\lesssim}


\define\eps{\varepsilon}

\define\la{\lambda}             \define\La{\Lambda}

\define\fA{{\frak A}}

\define\fW{{\frak W}}

\define\fg{{\frak g}}


\define\fw{{\frak w}}

\define\fz{{\frak z}}

\define\bbR{{\Bbb R}}

\define\bbZ{{\Bbb Z}}

\define\cC{{\Cal C}}

\define\cE{{\Cal E}}
\define\cF{{\Cal F}}
\define\cG{{\Cal G}}
\define\cH{{\Cal H}}

\define\cJ{{\Cal J}}
\define\cK{{\Cal K}}

\define\cT{{\Cal T}}
\define\cU{{\Cal U}}




\document

\head {\bf 1. Introduction}
\endhead

Let $G$ be a finite-dimensional
{\it  step two}  nilpotent group which we may identify with its Lie algebra
$\frak g$ by the exponential map.
We  assume that $\fg$ splits as  a direct sum 
$\fg=\fw\oplus\fz$ so that
$$[\fw,\fw]\subset \fz, \quad [\fw,\fz]=\{0\},$$
and  that 
 $\dim(\fw)=d$, $\dim(\fz)=\rk$.

Throughout we shall make   the following

\demo{\bf Nondegeneracy Hypothesis} For every nonzero linear functional  $\omega\in\fz^*$ the bilinear form
$$
\cJ_\omega:
\aligned  \fw\times\fw&\to \bbR
\\
(X,Y)&\mapsto \omega([X,Y])
\endaligned
$$ 
is nondegenerate.
\enddemo
Note that the skew symmetry of $\cJ_\omega$ and the  nondegeneracy hypothesis imply
that  $d$ is even.

There is a natural dilation structure relative to $\fw$ and $\fz$, namely
 for  $X\in \fw$ and $U\in \fz$ we consider the dilations
$$\delta_t: (X,U)\mapsto (tX, t^2 U).$$ 
With the identification of the Lie algebra with the group  $\delta_t$ becomes  an 
automorphism of the group.

In exponential coordinates  $(x,u)$, $x\in\Bbb R^d$, 
 $u\in \Bbb R^\rk$, the group  multiplication
is given by
$$
(x,u)\cdot (y,v)= (x+y, u+v+x^tJ y)
\tag 1.1$$
where
$x^tJy=(x^tJ_1y,\dots, x^t J_{\rk} y)\in\bbR^{\rk}$ and the
$J_i$ are skew-symmetric  matrices  acting on $\bbR^{d}$ (i.e. $J_i^t=-J_i$).
For $u\in \Bbb R^{\rk}$ we also form the skew-symmetric matrices
$J_u= \sum_{i=1}^{\rk} u_i J_i$ and the nondegeneracy hypothesis is 
equivalent with the invertibility of $J_u$ for all $u\neq 0$.

The most prominent examples are 
the  Heisenberg groups $H^n$ which arise when
$d=2n$, $\rk=1$ and $J=J_1$ is the standard symplectic matrix on
$\Bbb R^{2n}$. These belong to the class of {\it Heisenberg-type} groups
(termed  $H$-type groups in \cite {9}),
for which  $J_u^2 =-4|u|^2  I$, so that  the nondegeneracy hypothesis is 
clearly satisfied in this case.
Note that in general  $\rk$ has to  be small compared to $d$ 
(see \cite{9} where the connection with Radon-Hurwitz numbers is 
pointed out).
The class considered here  has been 
introduced by
M\'etivier  \cite{10} in his study of analytic hypoellipticity; the 
nondegeneracy assumption is termed ``Condition (H)'' in \cite{10}.
There are many groups which satisfy the nondegeneracy condition 
but which are not isomorphic
to a Heisenberg-type group; we give an example in \S7.

Let $\Sigma$ be a smooth  convex hypersurface in $\fw$ and let $\mu$ be a 
compactly
supported  smooth density 
on $\Sigma$. We make
the following
\demo{\bf Curvature Hypothesis} The Gaussian curvature of $\Sigma$ does not vanish on the
support of $\mu$.
\enddemo

Define the dilate $\mu_t$ by
$$\align
\inn{\mu_t}{f}
&=\int f(tx,0) d\mu(x).
\tag 1.2
\endalign
$$
We recall the definition of convolution $$
\align f*g(x,u)&=\int f(y,v) g((y,v)^{-1}\cdot (x,u)) dydv
\\
&=
\int f(y,v) g(x-y,u-v+x^tJy) dy dv
\tag 1.3
\endalign
$$
and define  for Schwartz-functions 
the maximal operator $M$ by
$$ Mf (x,u)=\sup_{t>0}|f*\mu_t (x,u)|.$$

We prove the following sharp result.

\proclaim{Theorem} Suppose $d>2$.
Then $M$ extends to a bounded operator on $L^p(G)$ if and only if 
$p>d/(d-1)$.
\endproclaim

\remark{Remarks}
(i) Other more ``regular''
 spherical maximal functions on the Heisenberg group have
  been considered  in \cite{2}, \cite{15}. In these papers the maximal
functions are generated by measures on hypersurfaces and the averaging operators are Fourier integral operators associated  to local
canonical graphs. In our work the maximal functions are generated
by measures on surfaces of codimension $\rk+1$, and the associated canonical relations project with fold singularities.

(ii)
A previous  result is due to
Nevo and Thangavelu \cite{12} who considered  the
case of spherical means on the noncentral part
of the  Heisenberg groups  ($\rk=1$) and obtained  $L^p$ boundedness in the smaller range
$p>(d-1)/(d-2)$, $d>2$.

(iii) Our theorem is an analogue  of Stein's theorem \cite{16} in the Euclidean case.
The necessity of the condition $p>d/(d-1)$ follows from 
the example in
\cite{16}; one tests $M$ on the function given by
$f(y,v)=|y|^{1-d}(\log|y|)^{-1}\chi(y,v)$ with a
 suitable cutoff function $\chi$. The $L^2$ methods in this paper  are not sufficient to establish $L^p$ boundedness for $p>2$ for the case
$d=2$
(that is, for  an extension of  Bourgain's result
 \cite{1} in the Euclidean case).

(iv) The result should remain true  for any
nilpotent Lie group of step $\le 2$;
{\it i.e.} the nondegeneracy hypothesis should not be necessary.
This is currently an open problem.

(v) As a corollary of the $L^p$ estimate for the maximal operator
one obtains the pointwise 
convergence result  $\lim_{t\to 0}\mu_t*f(x)= cf(x)$ 
almost everywhere, if $f\in L^p$ and $c=\int d\mu$.  Moreover the $L^p$ 
bounds of the maximal operator 
 are relevant for certain results  in ergodic theory,
 where one needs to have pointwise control for large $t$.

(vi) We use in an essential way the invariance of the subspace $\fw$ under the dilation group $\{\delta_t\}$.
Namely this implies a favorable bound for the principal  symbol of  $(d/dt)\mu_t$ on the fold surface 
of the associated canonical relation. A similar phenomenon was observed in \cite{11} for averages along light rays.

(vii) One can replace the measure on $\fw$ by a measure supported on a perturbed subspace  $\fW$ which is transversal to the center
but  no longer invariant under $\{\delta_t\}$;
then the phenomenon in the last remark does not occur. In the above coordinates $\fW$ is 
given as 
$$\fW=\{(x,\La x), x\in  \bbR^d\},
\tag 1.4$$
where $\La=(\La_{ij})$ is a $m\times d$ matrix. Define a measure $\mu^\Lambda_t$ by
$$
\inn{\mu^\La_t}{f}
=\int f(tx,t^2\Lambda x ) d\mu(x),
$$
we also set $\mu^\La:=\mu^\La_1$.
Consider the maximal operator $M^\La$ defined by 
$$M^\La f=\sup_{t>0}|f*\mu^\La_t|.
\tag 1.5
$$
 For general $\La$ we then prove the partial result
that $M^\La$ is bounded for $p>(3d-1)/(3d-4)$. We conjecture that  boundedness holds for $p>d/(d-1)$ which
 by our theorem holds true  for $\La=0$.

\endremark

{\it Notation:}
Given two quantities $A$ and $B$ we write
$A\lc B$  if there is a positive constant $C$, such that $A\le CB$.

\head{  \bf 2.
Preliminary decompositions}
\endhead
We shall present the argument  for the maximal operator $M^\La$ in (1.5).
We shall denote by $\La_j$ the $j^{\text{th}}$ column of $\Lambda$
and  by
 $\|\La\|$ the matrix norm of
 $\La$ with respect to the Euclidean norms on $\Bbb R^d$ and $\Bbb R^m$.
In what follows we shall always assume that $\|\La\|\le C_1$ for some fixed $C_1$ (and various bounds may depend on $C_1$). If $\|\La\|$ occurs explicitly in an estimate then we are interested in the behavior for $\La\to 0$,
as the case of our Theorem  corresponds to $\Lambda=0$.

We note that by    localizations and rotations in $\Bbb R^d$ 
  one can assume
 that $\mu$ has small support and that the projection of
$\Sigma$ to $\fw$  is given as  a graph $x_{d}=\Ga(x')$, $x'=(x_1,\dots,x_{d-1})$,
so that
$\nabla_{x'} \Gamma(0)=0$ 
and so that $\mu $ is supported in a small 
neighborhood
of $(0, \Gamma(0))$ (we may assume that  $|\nabla_{x'} \Gamma(x')|\le C_0^{-1}
c_0/100$) where
$c_0$, $C_0$ are defined in (5.10) below).
Note that a rotation has the effect  of replacing the matrices 
$J_i$ in the group law by
$Q^tJ_iQ$ with $Q\in SO(d)$.
We thus will
 need to prove an estimate which is uniform in these rotations.


Using the Fourier inversion formula for Dirac measures we may write
$$\mu^\La(x,u)= \chi(x,u)
 \iint 
e^{i\big(\sigma(x_d-\Gamma(x'))+\tau \cdot (u-\La x)\big)} 
 d\sigma d\tau
$$
where $\chi$ is a smooth compactly supported  function and 
the integral converges in the sense of oscillatory integrals (thus in the
sense of distributions).

 We split the integrals by introducing dyadic
 decompositions in
$(\sigma, \tau)$ and
then also in $\sigma$, when $|\sigma|<|\tau|$.

Let $\zeta_0\in C^\infty_0(\bbR)$ be an even function  so that
$\zeta_0(s)=1$ if $|s|\le 1/2$ and $\supp (\zeta_0)\subset(-1,1)$.
Also define $\zeta_1(s)=\zeta_0(s/2)-\zeta_1(s)$ and for 
$k\ge 1$, $1\le l<k/3$,
$$
\align
\beta_{0}(\sigma,\tau)&=\zeta_0(\sqrt{\sigma^2+|\tau|^2})
\tag 2.1.1
\\
\beta_{k,0}(\sigma,\tau)&=\zeta_1(2^{-k}\sqrt{\sigma^2+|\tau|^2})(1-\zeta_0(2^{-k}\sigma))
\tag 2.1.2\\
\beta_{k,l}(\sigma,\tau)&=\zeta_1(2^{-k}\sqrt{\sigma^2+|\tau|^2})\zeta_1(2^{l-k}\sigma)
\\
\widetilde \beta_{k}(\sigma,\tau)&=\zeta_1(2^{-k}\sqrt{\sigma^2+|\tau|^2})
\zeta_0(2^{[k/3]-k-1}\sigma).
\tag 2.1.3
\endalign
$$
Then observe that $$\beta_0+\sum_{k\ge 1}\big(\beta_{k,0}+\sum_{1\le l<k/3}\beta_{k,l}
+\widetilde \beta_{k}\big)=1,$$
and  for $k>0$ the function $\beta_{k,0} $ is supported
where $\sigma \approx 2^k$ and $|\tau|\lc 2^k$,
$\beta_{k,l} $ is supported
where $|\tau|\approx 2^k$ and $|\sigma|\approx 2^{k-l}$ and
$\widetilde \beta_k$ is supported where
$|\tau|\approx 2^k$ and $|\sigma|\lc  2^{2k/3}$.

Define
$$
\align
K^{0}(x,u)&= \chi(x,u)
\iint 
e^{i\big(\sigma(x_d-\Gamma(x'))+\tau \cdot (u-\La x)\big)} 
\beta_0(\sigma,\tau)d\sigma d\tau,
\tag 2.2.1
\\
K^{k,l}(x,u)
&= \chi(x,u)
\iint e^{i\big(\sigma(x_d-\Gamma(x'))+\tau \cdot (u-\La x)\big)} 
\beta_{k,l}(\sigma,\tau)d\sigma d\tau,\quad\text{ \   $0\le l<k/3$, }
\tag 2.2.2
\\
\widetilde K^{k}(x,u)&= \chi(x,u)
\iint 
e^{i\big(\sigma(x_d-\Gamma(x'))+\tau \cdot (u-\La x)\big)}  
\widetilde \beta_k(\sigma,\tau)d\sigma d\tau;
\tag 2.2.3
\endalign
$$
moreover  for $t>0$ define the  dilates
$$
[K^0_{t},
K^{k,l}_{t},
\widetilde K^{k}_{t}]
(x,u)= t^{-(d+2\rk)}
[K^0,
K^{k,l},
\widetilde K^{k}]
(t^{-1}x, t^{-2} u).$$
Note that $\mu^\La_t=
K^0_t+\sum_{k\ge 1}\big(K^{k,0}_t+\sum_{1\le l<k/3} K^{k,l}_t
+ \widetilde K^{k}_t\big)$.

Since
$K^0$ is a bounded compactly supported function
 the associated maximal function is controlled by the appropriate variant
of the Hardy-Littlewood maximal function  and therefore (\cite{17}) we have the inequality
$$
\big\|\sup_t |f* K^0_{t}|\big\|_p\le C_p \|f\|_p
$$
for $1<p\le \infty$.

Using known estimates for oscillatory integral operators
 with fold singularities 
and additional almost orthogonality estimates
we shall
derive in \S3 and  \S5  the following $L^2$  estimates.

\proclaim{Proposition 2.1}
Suppose $k>0$. Then
 for $0\le l<k/3$
$$\big\|\sup_{t}|f* K^{k,l}_{t}|\big\|_2 
\lc \sqrt{k} 2^{-k(d-2)/2}(1+\|\La\|2^l)^{1/2}\|f\|_2;
\tag 2.3
$$
moreover
$$
\big\|\sup_{t}|f* \widetilde K^{k}_t|\big\|_2 
\lc \sqrt{k} 2^{-k(d-2)/2} (1+\|\La\|2^{k/3})^{1/2}\|f\|_2
\tag 2.4
$$
\endproclaim

To obtain $L^p$ results we shall interpolate with weak type 
inequalities proved in \S6.

\proclaim{Lemma 2.2} Let $k>0$. For all $\alpha>0$
we have
$$
\meas\big(\{(x,u): \sup_{t>0} |f* K^{k,l}_t(x,u)|>\alpha\}\big)
\lc k 2^{k-l}(1+\|\La\|2^l) \alpha^{-1}\|f\|_1
\tag 2.5
$$
for $0\le l<k/3$ and
$$
\meas\big(\{(x,u): \sup_{t>0} |f* \widetilde K^{k}_t(x,u)|>\alpha\}\big)
\lc k 2^{2k/3}(1+\|\La\|2^{k/3}) \alpha^{-1}\|f\|_1.
\tag 2.6
$$
\endproclaim

We interpolate by the real method and obtain
\proclaim{Corollary 2.3}
Suppose $1<p\le 2$ and $k>0$. Then
 for $0\le l<k/3$
$$\big\|\sup_{t}|f* K^{k,l}_{t}|\big\|_p 
\le C_p {k}^{1/p} 2^{-k(d-1-d/p)}2^{-l(2/p-1)}(1+\|\La\|2^l)^{1/p}\|f\|_p;
\tag 2.7
$$
moreover
$$
\big\|\sup_{t}|f* \widetilde K^{k}_t|\big\|_p
\le C_p k^{1/p}
2^{-k(d-4/3-d/p+2/3p)}
 (1+\|\La\|2^{k/3})^{1/p}\|f\|_2.
\tag 2.8
$$
\endproclaim

Now if $p<2$ we may sum   in $k$ and $l$ and see  that  $M^\La$ is $L^p$ bounded if 
$d-4/3-d/p+1/(3p)>0$ which is equivalent to  $p>(3d-1)/(3d-4)$
(showing the estimate mentioned in remark (vii) in the introduction).
If $\Lambda=0$ we get a better bound, namely that
  $L^p$ boundedness holds if 
$d-1-d/p>0$ or  $p>d/(d-1)$. This proves our main Theorem.

\head{\bf 3. Square functions and almost orthogonality
}
\endhead

It is advantageous to introduce cancellation in the above kernels,
modulo small  acceptable errors.
Indeed 
$$
\Big|\iint
K^{k,l}(x,u) dx du\Big|+\Big| \iint \widetilde K^{k}(x,u) dx du\Big|
\le C_N 2^{-kN},
$$
for all $N=0,1,\dots$,
and this estimate follows by an integration by parts in the 
$(x,u)$ variables.
Thus there is a $C^\infty_0$ function $b$ which is  equal to $1$  on
$\supp (\chi)$, and constants $\ga_{k,l}$, $\gamma_k$ so that
$$
\aligned
&\iint
K^{k,l}(x,u) dx du = \ga_{k,l}  \iint b(x,u) dx du
\\
&\iint
\widetilde K^{k}(x,u) dx du = \ga_{k}  \iint b(x,u) dx du
\endaligned
\tag 3.1
$$
where
$$|\ga_{k}|+|\ga_{k,l}|
\le C_N 2^{-kN}.
\tag 3.2
$$
We define
$$
\align \cK^{k,l}(x,u)
& = K^{k,l}(x,u) -\ga_{k,l}  b(x,u)
\tag 3.3.1
\\
\widetilde \cK^{k}(x,u)&=
\widetilde K^{k}(x,u) - \ga_{k} b(x,u)
\tag 3.3.2
\endalign
$$
and denote by
$\cK^{k,l}_{t}$, $\widetilde \cK^{k}_t$ their dilates, as before.
Then the functions
$\cK^{k,l}_{t}$, $\cK^{k}_t$  have integral zero.

Since the maximal operator generated
by the kernel $b$ (with nonisotropic dilations) is bounded
by the nonisotropic Hardy-Littlewood maximal operator we see that
for $1<p\le \infty$
$$
\big\|\sup_t |
f*
(\cK^{k,l}_{t}-
K^{k,l}_{t})|\big\|_p\le C_{N,p} 2^{-kN} \|f\|_p.
$$

Now in order to deal with the main term we shall use the following
standard lemma in the subject  which is an immediate consequence 
of a similar one stated in  \cite{17, p.499}.

\proclaim{Lemma 3.1}
Suppose that
$$\align &\sup_{s\in [1,2]}\Big(\sum_{n\in \bbZ}
\big\|F_n(\cdot,s)\big\|_2^2\Big)^{1/2} \le A_1
\\
&\sup_{s\in [1,2]}\Big(\sum_{n\in \bbZ}
\big\|\frac{\partial F_n}{\partial s} (\cdot,s)\big\|_2^2\Big)^{1/2}
\le A_2.
\endalign
$$
Then
$$
\Big\|
\sup_n\sup_{s\in [1,2]} |F_n(\cdot,s)|
\Big\|_2\le C( A_1+ \sqrt{A_1A_2}).
$$
\endproclaim

We omit the proof.
Using Lemma 3.1 one sees that the estimates
$$\align
\big\|\sup_{t}|f* \cK^{k,l}_{t}|\big\|_2
&\lc \sqrt{k} 2^{-k(d-2)/2}
(1+\|\La\|2^l)^{1/2}
\|f\|_2\\
\big\|\sup_{t}|f* \widetilde K^k_{t}|\big\|_2 &\lc
\sqrt{k} 2^{-k(d-2)/2}
(1+\|\La\|2^{k/3})^{1/2}
\|f\|_2
\endalign
$$
follow from the following estimates which are uniform in $s\in [1,2]$.
$$
\align
&\Big(\sum_n\big\|f* \cK^{k,l}_{2^n s}\big\|_2^2\Big )^{1/2} \lc \sqrt k 2^{-k(d-1)/2} 2^{l/2}\|f\|_2
\tag 3.4
\\
&\Big(\sum_n
\Big\|f* \Bigl[t \frac{\partial}{\partial t}\cK^{k,l}_{t}\Bigr]_{t=2^ns}
\Big\|_2^2\Big )^{1/2} 
\lc 
 \sqrt k
2^{-k(d-3)/2} 2^{-l/2}
(1+\|\La\|2^{l})
\|f\|_2, 
\tag 3.5
\endalign
$$
for $l<k/3$,  and
$$
\align
&\Big(\sum_n\big \|f*\widetilde  \cK^{k}_{2^ns}\big\|_2^2\Big )^{1/2} \lc
 \sqrt k 2^{-k(d-1)/2+k/6} \|f\|_2
\tag 3.6
\\
&\Big(\sum_n
\Big\|f*
\Bigl[t\frac{\partial}{\partial t}\widetilde \cK^{k}_{t}
\Bigr]_{t=2^ns}
\Big\|_2^2\Big )^{1/2} 
\lc \sqrt k
2^{-k(d-3)/2-k/6}
(1+\|\La\|2^{k/3})
\|f\|_2. \tag 3.7
\endalign
$$

Note by scaling that it suffices to prove these estimates for $s=1$.
We shall first
use the cancellation of the kernels $\cK^{k,l}_{2^n s}$ and
 $\widetilde \cK^{k}_{2^n s}$
to show certain
almost orthogonality properties (for the sums in $n$) and then
 we use stronger estimates for oscillatory integrals to establish
decay estimates for fixed $n$.

\subheading {An almost orthogonality lemma}
We first state a simple and  presumably well known
consequence
of the Cotlar-Stein Lemma.

\proclaim{Lemma 3.2}
 Suppose $ 0<\eps<1$, $A\le B/2$ and let
 $\{T_n\}_{n=1}^\infty$ be a sequence  of bounded operators on a
Hilbert space  $H$ so that the operator norms satisfy
$$
\|T_n\|\le A
\tag 3.8
$$
and
$$\|T_n T_{n'}^*\|
\le B^2 2^{-\eps|n-n'|}.
\tag 3.9$$

Then for all $f\in H$
$$
\Big(\sum_{n=1}^\infty\|T_n f\|^2\Big)^{1/2} \le  C A
\sqrt{\eps^{-1}\log(B/A)} \|f\|.
\tag 3.10
$$
\endproclaim

\demo{\bf Proof} For $N\ge 1$ consider the operator 
$$\cT_N: H\to \ell^2(H)$$  which maps $f$ to the sequence
$(T_1 f,\dots, T_N f,0,0,\dots)$.
Now $\|\cT_N\|= \|\cT_N^* \cT_N\|^{1/2}$ where
$\cT_N^*\cT_N:H\to H$
is given by
$$
\cT_N^*\cT_N f = \sum_{n=1}^N T_n^* T_n f.
$$
We let $S_n =T_n^* T_n$ and observe that
$$\align&
\|S_{k}^* S_l\|=\|S_{k} S_l^*\|
=
\|T_k^*T_kT_l^*T_l\|
\\&
\le \|T_k^*\|\|T_k T_l^*\|\|T_l\|
\le A^2\min \{ A^2, B^22^{-|k-l|\eps}\}.
\endalign
$$
The standard Cotlar-Stein Lemma \cite{17} gives
$$\|\cT_N^*\cT_N\|\le \sum_{m=\infty}^\infty
\max\big\{ \sup_{k-l=m}\|S_k^* S_l\|^{1/2}, \sup_{k-l=m}\|S_k S_l^*\|^{1/2}\big\}
$$
and thus
$$\align\|\cT_N\|^2 &\le A\sum_{m=-\infty}^\infty \min\{ A, B2^{-|m|\eps}\}
\\&\le  C^2\eps^{-1} A^2\log(B/A).\endalign
$$
Thus $\|\cT_N f\|_{\ell^2(H)}$ is dominated by the right hand side of 
(3.10), and the assertion follows by taking the limit as $N\to \infty$.
\qed\enddemo

\remark{Remark} We proved  Lemma 3.2 by using the statement of the 
Cotlar-Stein Lemma. Using the proof
of the Cotlar-Stein Lemma one can also show the 
following more general fact:
If
$\|T_n T_{n'}^*\|\le \alpha^2(n-n')$ then
$$\Big(\sum_{n=1}^N\|T_n f\|^2\Big)^{1/2} \lc
\Big(\sum_{j\in \bbZ}|\alpha(j)|^2\Big)^{1/2}\|f\|.$$
Of course, Lemma 3.2 is an immediate consequence of this inequality.\endremark

\subheading{Almost orthogonality estimates}
Here we wish to apply Lemma 3.2 to convolutions on groups.
If $Tf=f*g$ we first note that  its adjoint is given by $T^* f=f*g^*$ where
$g^*= \overline{g(\cdot^{-1})}$. 
Moreover using Minkowski's inequality and the  unimodularity of nilpotent Lie groups 
one obtains the standard convolution inequality 
$$\|f*g\|_2\le \|g^*\|_1\|f\|_2=\|g\|_1\|f\|_2.$$

We now fix $k,l$ and $s\in [1,2]$  and derive almost orthogonality properties for the operators of convolution with $\cK^{k,l}_{2^n s}$.

Notice that for $n\le 0$ the function
$\cK^{k,l}_{2^n s}$
is supported in a (small) ball
of radius $C 2^{n}$ (in fact in a smaller nonisotropic ball).
Moreover we have
$|\nabla_{y,v} \cK^{k,l}_{s}(y,v)|\le 2^{k(\rk+2)}$ and using the
cancellation of $\cK^{k,l}_{2^n s}$ we obtain
$$
|\cK^{k,l}_{s}* (\cK^{k,l}_{2^n s})^*(x,u)|
\lc 2^{k(\rk+2)} 2^{n} \text{ if $n\le 0$}.
$$
By scaling and applying Schur's Lemma we obtain
$$
\big\|f*
\cK^{k,l}_{2^{n'}s}* (\cK^{k,l}_{2^n s})^*\big\|_2\lc
2^{k(\rk+2)} 2^{-|n-n'|}\|f\|_2\tag  3.11
$$
first for $n\le n'$ and then by taking adjoints also for $n< n'$.
This and the following estimates are uniform in $s\in [1,2]$.

 Similarly we get
$$\Big\|f*s
\frac{\partial \cK^{k,l}_{2^{n'}s }}{\partial s}
* s\frac{\partial ({\cK^{k,l}_{2^{n'}s}})^*}{\partial s}\Big\|_2\lc
2^{k(\rk+4)} 2^{-|n-n'|}\|f\|_2
\tag 3.12
$$
and also
$$
\align
&\big\|f*
\widetilde \cK^{k}_{2^{n'}s}
* (\widetilde \cK^{k}_{2^ns})^*\big\|_2\lc 2^{k(\rk+2)} 2^{-|n-n'|}\|f\|_2.
\tag 3.13
\\
&\Big\|f*s \frac{\partial\widetilde  \cK^{k}_{2^n s}}{\partial s}*
s\frac{\partial(\widetilde  \cK^{k}_{2^{n'}s})^*}{\partial s}\Big \|_2
\lc 2^{k(\rk+4)}
2^{-|n-n'|}\|f\|_2.
\tag 3.14
\endalign
$$
In \S5 we shall prove the inequalities
$$
\align
&\|f*K^{k,l}
 \|_2\lc 2^{-k (d-1)/2}2^{l/2} \|f\|_2 
\tag 3.15
\\&\Big\|f*\Bigl[\frac{\partial K^{k,l}_s}{\partial s}\Bigr]_{s=1}
\Big \|_2\lc
2^{-k(d-3)/2} 2^{-l/2}(1+\|\La\|2^l)  \|f\|_2
\tag 3.16
\endalign
$$
for $l<k/3$, and
$$
\align
&\|f*\widetilde K^{k}
 \|_2\lc 2^{-k (d-1)/2}2^{k/6} \|f\|_2 
\tag 3.17
\\&\Big \|f*\frac{\partial \widetilde K^{k}_s}{\partial s}
\big|_{s=1}
\Big \|_2\lc
2^{-k(d-3)/2} 2^{-k/6}
(1+\|\La\|2^{k/3})\|f\|_2.
\tag 3.18
\endalign
$$
By scaling and by (3.2)
 the same inequalities hold with
$K^{k,l}$ and $\widetilde  K^k$ replaced by
$\cK^{k,l}_t$ and  $\widetilde \cK^k_t$ and with
$\partial_s K^{k,l}$,
$\partial_s \widetilde  K^{k}$ replaced by
$\partial_s  \cK^{k,l}_{2^ns }$,
${\partial_s\widetilde \cK^{k}_{2^n s}}$, for $1\le s\le 2$.

Now the inequality (3.4)  follows from (3.15) and (3.11) 
 if we apply Lemma  3.2 with $A= 2^{-k(d-1)/2} 2^{l/2}$ and $B=2^{k(m+4)}$.
Similarly (3.5) follows from (3.16) and (3.12), (3.6) from
(3.17) and (3.13), and (3.7) from (3.18) and (3.14).

The next two sections are concerned with 
the derivation of inequalities (3.15-18).
 

\head{\bf 4. Preliminaries   on oscillatory integral operators with
folding canonical relations}\endhead
We shall reduce matters to estimates for 
oscillatory integral operators whose canonical relations  have two-sided fold
singularities. We consider localizations near the fold surface and the 
 estimate goes back to Phong and Stein
\cite{13} for certain conormal operators in the
plane;  the general case is implicit in  Cuccagna's paper \cite{3}.
For the version needed here we refer to  \cite{6}.

Let $\Omega\in \Bbb R^n\times\Bbb R^n$ be an open set and let $\Gamma$ be an open set in some finite dimensional space.
We consider phases $\varphi(x,y,\gamma)$ and
amplitudes
$a_\lambda(x,y,\gamma)$, $(x,y,\gamma)\in
\Omega\times\Omega\times \Gamma$, and assume that
$$
\align
&|\partial_x^\alpha \partial_y^\beta \varphi(x,y,\gamma)|\le C
\tag 4.1
\\
&|\partial_x^\alpha \partial_y^\beta a_\lambda(x,y,\gamma)|\le C
\lambda^{(|\alpha|+|\beta|)/3}
\tag 4.2
\endalign
$$
say, for all multiindices $\alpha, \beta$ with $|\alpha|, |\beta|\le 10 n$,
with uniform bounds in $\Omega\times \Gamma$; we also assume that all derivatives depend continuously on the parameter $\gamma$.

We shall assume that 
$$\cC_\varphi=\{(x,\varphi_x, y,-\varphi_y)\}$$
 is a folding canonical
relation, i.e. for each point $P_0=(x_0,y_0,\gamma_0)$ we have
$$
\rank \varphi_{xy}''(P_0)\ge n-1,
\tag 4.3
$$
and
for unit vectors $U$, $V$
$$
\align
\varphi_{xy}''(P_0)V=0 \quad  &\implies \quad \big|\inn V{\nabla_y} \det \varphi_{xy}''\big|\ge c,
\tag 4.4
\\
U^t\varphi_{xy}''(P_0)=0 \quad  &\implies \quad \big|
\inn U{\nabla_x} \det \varphi_{xy}''\big|\ge c,
\tag 4.5
\endalign
$$
for some $c>0$.

We consider the oscillatory integral operator
$T_{\la}[b] $ defined by
$$
T_{\la}[b] f(x)=\int e^{i\lambda \varphi(x,y,\gamma)} b(x,y,\gamma) f(y) dy
$$
which is bounded on all $L^p$ if $b$ is bounded and  compactly supported.
We shall take for $b$ certain localizations of the symbol in terms of the
size of $\det\varphi_{xy}''$.
Let $\eta$ be smooth and compactly supported in $(-1,1)$ so that
$\eta(s)=1$ for $|s|\le 1/2$ and set
$$\beta_l(x,y,\gamma)=\eta(2^l\det\varphi_{xy}''(x,y,\gamma))
-\eta(2^{l+1}\det\varphi_{xy}''(x,y,\gamma)),
$$ so that $\beta_l $ localizes to the set where
$|\det\varphi_{xy}''|\approx 2^{-l}$. We also define
$$\zeta_\la(x,y)=1-\sum_{2^l<\la^{1/3}} \beta_l(x,y)$$
so that $|\det \varphi_{xy}''|\lc \la^{-1/3}$ on $\supp(\zeta_\la)$.

Then
there is a neighborhood
$\cU$ of $(x_0,y_0,\gamma_0)$ so that for all $a_\lambda$ satisfying 
(4.2), supported in $\cU$ the following 
estimates hold for the operator norms:

$$
\big\| T_\lambda[ a_\la \beta_l]\big\|_{L^2\to L^2} \le C_1
2^{l/2} \lambda^{-n/2},
\qquad 2^l\le \la^{1/3}\tag 4.6
$$
and
$$
\big\| T_\lambda[ a_\la \zeta_\la]\big\|_{L^2\to L^2} \le C_1
\lambda^{1/6-n/2}.\tag 4.7
$$
These estimates are a consequence  of Theorem 2.1 in \cite{6}.

\head{\bf 5. Reduction to oscillatory integral operators}
\endhead
We now consider  the operator of convolution with  
$K^{k,l}$ and give the 
proof of  the bound (3.15). The operator $\partial_s K^{k,l}$ is more 
singular, but its estimation is rather analogous,
so we shall point out the modifications needed for (3.16)
 at the end of this section. The estimations for $\widetilde K^k$ and $\partial_s 
\widetilde K^k_s$ will be similar.

Since
$K^{k,l}$ is compactly supported in a fixed neighborhood
we may use the  translation invariance to reduce to the case that
 $f$ is also compactly supported in a
fixed neighborhood of the origin. Thus it suffices 
to show the desired bound for the operator
with Schwartz kernel
$$
\chi_1(x,u) K^{k,l}(x-y, u-v+x^tJy) \chi_2(y,v),
\tag 5.1
$$
for suitable compactly supported smooth functions $\chi_1$ and $\chi_2$.
In what follows we set $\lambda =2^k$ and then by a change of variables
the kernel (5.1) can be written as
$$
H^{\lambda,l}(x,u,y,v)= \lambda^{m+1}\iint 
e^{i\la\phi(x,u,y,v,\sigma,\tau)}
\chi_0(x,u,y,v) \eta_l( \sigma,\tau) d\sigma d\tau
\tag 5.2
$$
where
$$\phi(x,u,y,v, \sigma,\tau)=\sigma(x_d-y_d-\Ga(x'-y'))+
\tau\cdot(u-v+x^tJy-\La(x-y))$$ 
and 
where $|\tau|\approx 1$ and $|\sigma|\approx 2^{-l}$ on the
support of $\eta_l$; specifically 
$$\eta_l(\sigma,\tau)=\zeta_1(\sqrt{\sigma^2+|\tau|^2})\zeta_1(2^l\sigma),
$$ and
$
\chi_0(x,u,y,v) = \chi_1(x,u) \chi(x-y, u-v+x^tJy)\chi_2(y,v)
$.

\remark{Notation} We let  $P: \Bbb R^d \to \bbR^{d-1}$ be the
linear  map with $Pe_i=e_i$, $i=1,\dots,d-1$ and $Pe_d=0$. We also
 use the notation $P$ for the $(d-1)\times d$ matrix 
$$P=\pmatrix I&0\endpmatrix$$ 
and $P^t$ for its transpose. 
\endremark

\subheading{Stationary phase calculations}
We wish to apply stationary phase arguments to reduce matters to the estimation of an
oscillatory integral operators without frequency variables (see {\it e.g.} 
 the general discussion in \cite{5}).

We shall apply a 
scaled Fourier transform on $\Bbb R^{m+1}$,
in the $(x_d,u)$ variables. 
Define
$$\Cal F_\lambda g(x',x_d,u)=
\iint e^{-i\lambda (x_d z_d+u\cdot w)} g(x',z_d,w) dz_d dw;
$$
then $(\lambda/2\pi)^{(m+1)/2}\cF_\la$ is a unitary operator and thus,
if $\Cal H^{\la,l}$ denotes the operator with Schwartz kernel 
$H^{\la,l}$  we have to prove
that $\cF_\la \cH^{\la,l}$ maps $L^2$ to itself with operator norm
$O(\la^{-(d+m)/2} 2^{l/2})$.
Let $\chi_3(x_d,u)$ denote a
smooth compactly supported function which is equal to one whenever $|x_d|+|u|\le 10$, and
define $\Cal F_{\lambda,1}$ by
$$\Cal F_{\lambda,1}
 g(x',x_d,u)=\chi_3(x_d,u)
\iint e^{-i\lambda (x_d z_d+u\cdot w)}
 g(x',z_d,w)
dz_d dw;
$$
moreover  let
$\Cal F_{\lambda,2}=\Cal F_{\lambda}-\Cal F_{\lambda,1}$.
Then the Schwartz kernel  of
$\Cal F_{\lambda,1} H^{\la,l}$
is given by
$$
\la^{m+1}
\int e^{i\la \Psi(x,u, y,v,\theta)}b_{l}(x,u,y,v, \theta) d\theta
\tag 5.3
$$
where  with $$\theta=(z_d,w,\sigma,\tau)$$  the phase function $\Psi$
 is given by
$$\align \Psi(x,u,y,v,\theta)=&-x_dz_d-u\cdot w+
\sigma\big(z_d-y_d-\Gamma(x'-y')\big)
\\&+\tau^t \big(w-v+
\La P^t(x'-y')+\La_d( z_d-y_d)+
({x'}^t,z_d)Jy\big),
\endalign
$$
and the amplitude is given by
$$b_l(x,u,y,v,\theta)=\chi_3(x_d,u)
\chi_0(x',z_d, y,w)
\eta_l(\sigma,\tau).
$$

For the error term  $\Cal F_{\lambda,2} H^{\la,l}$
we have a similar formula, only with $\chi_3 $ replaced by $1-\chi_3$.
Then in view of the support properties of $(1-\chi_3)$ we see that
$|\nabla_{z_d,w}\Psi|\ge |x_d|+|u|$ on $\supp(1-\chi_3)$ and by integration by parts with respect to the $(z_d,w)$ variables we see that the kernel of
$\Cal F_{\lambda,2} H^{\la,l}$ is bounded by
$C_N\la^{m+1-N} (|x_d|+|u|)^{-N}$. Moreover this kernel is supported on a set
where
$|x_d|+|u|\ge 1$ and
where $|x'|+|y|+|v|\le C$. Thus, with  an obvious application of
 Schur's Lemma  we conclude that the operator
$\Cal F_{\lambda,2} H^{\la,l}$ is bounded on $L^2$ with operator norm
$O(\lambda^{-N})$ for any $N$.

We return to the main term $\Cal F_{\lambda,1} H^{\la,l}$ 
and it remains to be  shown that
$$\|\Cal F_{\la,1} \cH^{\la,l}\|\lc 2^{l/2} \la^{-(d+m)/2}.
\tag 5.4
$$

 Note that for fixed $(x,u,y,v)$ the
phase function $\Psi$ is a polynomial of degree $\le 2$ in the 
$\theta$ variables and that the Hessian 
$\Psi_{\theta\theta}''$ is nondegenerate.

Indeed,
$$
\aligned
\Psi_{z_d}'&= -x_d+e_d^tJ_\tau y+\sigma+\tau^t \La_d
\\
\Psi_{w}'&=\tau-u
\\
\Psi_{\tau}'&= w-v+({x'}^t,z_d)Jy +\La P^t(x'-y')+\La_d(z_d-y_d)
\\
\Psi_\sigma'&=z_d-y_d-\Ga(x'-y')
\endaligned
\tag 5.5
$$
and with $\Xi$ denoting the column vector in $\Bbb R^m$  with coordinates $\Xi_i=e_d^t J_i y+\La_{id}$ we have
$$
\Psi_{\theta\theta}''=
\pmatrix
0&0&\Xi^t&1
\\
0&0&I&0
\\
\Xi&I&0&0
\\
1&0&0&0
\endpmatrix .
$$

Clearly the linear equations $\Psi_\theta=0$ have a unique solution
$\theta_{\text{crit}}=
(z_d,w,\tau,\sigma)_{\text{crit}}$, with
$$\align
(z_d)_{\text{crit}}(x,u,y,v)&=
y_d+\Ga(x'-y')
\\
(w_i)_{\text{crit}}(x,u,y,v)
&=v_i-({x'}^t, y_d+\Ga(x'-y'))J_i y-e_i^t\La P^t(x'-y')-
\La_{id}\Ga(x'-y')
\\
(\tau_i)_{\text{crit}}(x,u,y,v)&=u_i
\\
{\sigma}_{\text{crit}}(x,u,y,v)&= x_d- \sum_{i=1}^m u_i(
e_d^tJ_iy+ \La_{id})
\endalign
$$
 and we can apply the method of stationary phase (with respect to the $2(m+1)$ frequency
 variables $\theta$).
 Setting
$$\align
\Phi&(x,u,y,v):=\Psi(x,u,y,v,\theta_{\text{crit}}(x,u,y,v))
\\&=-x_d(y_d+\Gamma(x'-y'))-\sum_{i=1}^m
u_i\big(  v_i
-({x'}^{t},y_d+\Ga(x'-y'))J_iy 
-\La_{id}\Ga(x'-y')-e_i^t \La P^t(x'-y')
\big)
\tag 5.6
\endalign
$$
we obtain that
$$
\align
\la^{m+1}
\int e^{i\la \Psi(x,u,y,v,\theta)}b_{l}(x,u,y,v,\theta) d\theta&=
e^{i\la\Phi(x,u,y,v)}
\sum_{j=0}^{N-1}\cE_j^l(x,u,y,v)\lambda^{-j}
+ R_N^{\la,l}(x,u,y,v)
\tag 5.7
\endalign
$$
where
$$\multline
\cE_j^l(x,u,y,v)=\\ (2i)^{-j}
\big(\det (\Psi_{\theta\theta}(x,y,u,v,\theta_{\text crit}(x,u,y,v)
)/2\pi i
\big)^{-1/2}
\frac 1{j!}\inn{ \Psi_{\theta\theta}^{-1}D_\theta}{D_\theta}^j
b_l(x,u,y,v, \theta)\Big|_{\theta=\theta_{crit}(x,u,y,v)}
\endmultline
\tag 5.8
$$
and
$$
| R_N^{\la,l}(x,u,y,v)|\le C_N
\|b_l\|_{L^2_{m+2+2N}}  \la^{-N}\le C_N'
2^{l(m+2+2N)}\la^{-N}.
\tag  5.9
$$
Here we have applied Lemma 7.7.3 in \cite{7}.

Since $2^l\le \la^{1/3}$ the error term $R_N^{\la,l}$ 
(which is compactly supported)
 defines a bounded operator
on $L^p$
with norm $O(\la^{-(2m+1+N)/3})$  which for large $N$ is much
 better than the desired bound in (5.4).

\proclaim{Claim 5.1} The operators with kernels
$\la^{-j}\cE_j^l(x,u,y,v) e^{i\la\Phi(x,u,y,v)}$ have $L^2$ operator norm
$O(\lambda^{-(d+m)/2-j/3} 2^{l/2})$
\endproclaim

This clearly implies (5.4).

\subheading{Geometry of the canonical relation}

We consider the canonical relation
$\Cal C_\Phi=(x,u,\Phi_x, \Phi_u; y,v,-\Phi_y,-\Phi_v)$ and the singularities of the
maps $p_L:(y,v)\mapsto(\Phi_x,\Phi_u)$,
$p_R:(x,u)\mapsto(\Phi_y,\Phi_v)$.
It is our objective to check the analogues of (4.3-4.5) and we will have to verify a few 
elementary linear algebra facts.


Let $A$ denote the $(d-1)\times(d-1)$ matrix $\Gamma''(x'-y')$ and let $B$ denote the column
 vector $\Gamma'(x'-y')\in \Bbb R^{d-1}$; recall that we may assume that $|B|$ is small. Indeed if
$$
\align c_0&=\min_{u\in S^{m-1}} \|J_u^{-1}\|^{-1}
\tag 5.10.1
\\
C_0&=\max_{u\in S^{m-1}}\|J_u\| 
\tag 5.10.2
\endalign
$$
we may assume that $$\|B\|\le C_0^{-1}c_0/100.$$

Now $p_L$ is explicitly given by

$$
\aligned
\Phi_{x'}&= -x_d \Ga'(x'-y')
+ PJ_u y+ \Gamma'(x'-y') e_d^tJ_u y+ u^t\Lambda_d \Gamma'(x'-y')+
u^t\Lambda P^t
\\
\Phi_{x_d}&=-y_d-\Gamma(x'-y')
\\
\Phi_{u_i}&= -
\big(  v_i
-({x'}^{t},y_d+\Ga(x'-y'))J_iy
-e_i^t \La P^t(x'-y')-\La_{id}\Ga(x'-y')\big).
\endaligned
$$

We compute the  differential  $Dp_L$ as 
$$
\Phi_{(x,u),(y,v)}''
=\pmatrix
(x_d- e_d^tJ_u y-u^t \La_d)A +PJ_uP^t  +B e_d^tJ_uP^t 
& PJ_ue_d&0
\\
B^t&-1&0
\\
C& c &I
\endpmatrix
\tag 5.11
$$
where $I$ is an $m\times m$ identity matrix and
$C$ is $m\times (d-1)$ matrix with rows $C_i= {x'}^{t}PJ_iP^t +
y_d  e_d^tJ_i P^t -  (e_d^t J_iy+\La_{id}) B^t+e_i^t\La P^t 
 +\Ga(x'-y')  e_d^t J_iP^t $ and $c$ is the column in $\bbR^m$ 
with $c_i=({x'}^{t},0)J_ie_d
+e_d^tJ_i y$.
In this calculation the skew symmetry of the $J_i$ is used.

We now  compute  the determinant of (5.11)
and obtain
$$\det
\Phi_{(x,u),(y,v)}''=
(-1)^d \det
\Big((x_d- e_d^tJ_u y-u^t\La_d)A +PJ_uP^t  +
E(B)\Big)
\tag 5.12
$$
where $$E(B)=
B e_d^tJ_uP^t  +PJ_ue_dB^t.
\tag 5.13$$

Here we used the factorization
$$\pmatrix
\sigma A +PJ_uP^t  +B e_d^tJ_uP^t 
& PJ_ue_d
\\
B^t&-1\endpmatrix=
\pmatrix
\sigma A +PJ_uP^t  +E(B)
& PJ_ue_d
\\0&-1\endpmatrix\pmatrix I&0\\-B^t&1\endpmatrix.
$$

Note that  $E(B)$ is  a skew-symmetric $(d-1)\times(d-1)$ matrix
and so is
 $PJ_uP^t +E(B)$. Thus, since $d-1$ is odd, 
 the rank  of $PJ_uP^t +E(B)$ is at most $d-2$, and
the following lemma shows that for
 small $B$ the rank is equal to $d-2$.

\proclaim{Lemma 5.2}
Suppose that
$$\|B\|\le \frac{c_0}{4C_0}.$$
Then the following holds:

(i) If $W \in \ker (PJ_uP^t +E(B)) $
then
$$|e_d^t J_u P^t W|\ge \frac{c_0}2 \|W\|.
\tag 5.14 $$

(ii) $\dim \ker (PJ_uP^t +E(B))=1$.

(iii) If $X$ belongs to the orthogonal complement of
$\ker (PJ_uP^t +E(B)) $   then
$$
\|(PJ_uP^t +E(B)) X\|\ge \frac{c_0}2  \|X\| .
\tag 5.15
$$
\endproclaim

\demo{Proof} Observe that 
$$\|E(B)\|\le 2C_0\|B\|.$$
Thus if $W \in \ker (PJ_uP^t +E(B)) $  and $\|W\|=1$ then
$$
\align
1=
\|P^t W\|&\le \|J_u^{-1}\|\|J_u P^t W\|
\\ &\le
\|J_u^{-1}\|\big(|e_d^tJ_u P^t W|+
\|PJ_u P^t W\|\big)
=
\|J_u^{-1}\|\big(|e_d^tJ_u P^t W|+
\|E(B)W\|\big)
\\&\le
c_0^{-1}\big(|e_d^tJ_u P^t W|+
2C_0\|B\|\big)
\endalign
$$
and thus, if $\|B\|\le c_0/4C_0$ we obtain
$|e_d^tJ_u P^t W|\ge c_0/2$ which is (5.14).

Let $S_u=J_u +E(B)$. Since $S_u$ is skew symmetric, it can be diagonalized over $\Bbb C$, and the eigenvalues are imaginary. The bounds (5.10.1/2) are still valid if $J_u^{-1} $ is acting as a linear transformation on $\Bbb C^d$.
  Let $\eta\in 
\Bbb C^d$  be a  unit eigenvector of $S_u$ so that $S_u \eta=i\la \eta$ and 
$\|\eta\|=1$; then
$$
|\la|=\|S_u\eta\|\ge \|J_u\eta\|-\|E(B)\eta\|\ge c_0-\|E(B)\|
\ge c_0-2C_0\|B\|\ge \frac{c_0}2
$$ by assumption on $B$.
Hence $|\la|\ge c_0/2$ for every eigenvalue $i\la$ of $S_u$. In particular $S_u$ is nondegenerate. But then $PS_uP^t=PJ_uP^t+E(B)$ has 
rank $d-2$ and therefore a one-dimensional kernel and all nontrivial eigenvalues of $S_u$ are also eigenvalues of
 $PS_uP^t$. This implies for vectors $X$ orthogonal to the kernel of $PS_uP^t$ that
$$PS_uP^t X\ge \frac{c_0}2 \|X\|$$   which is (5.15).
\enddemo

\proclaim{Lemma 5.3} 
Let $\fA$ be a symmetric 
 positiv definite  matrix on $\Bbb R^{n}$ and let $S$
be a skew-symmetric matrix on $\Bbb R^{n}$.
Then:

(i)  For all $\sigma\neq 0$, the matrix $\sigma \fA+S$ is
invertible
and the inverse satisfies the bounds
$$
\|(\sigma \fA+S)^{-1}\|\le 
|\sigma|^{-1}\|\fA^{-1}\|.
\tag 5.16
$$

(ii)  If $S$ is invertible then $\sigma \fA+S$ is invertible for all $\sigma$ and we have the bound
$$
\|(\sigma \fA+S)^{-1}\|\le 
2\|S^{-1}\| \quad \text{ if } |\sigma|\le
\big(2\|\fA\| \|S^{-1}\|\big)^{-1}.
\tag 5.17
$$
\endproclaim
\demo{Proof}

For a unit vector $e$ in $\bbR^{n}$ we get
$$\align
\|(\sigma \fA+S)e\|&\ge|
\inn
{(\sigma \fA+S)e}{e}|=
|\inn
{\sigma \fA e}{e}|\ge |\sigma|\|\fA^{-1}\|^{-1}.
\endalign
$$
Here we have used that by the skew symmetry of $S$ we have
$\inn{Se}{e}=0$, and also that
  $\|\fA^{-1}\|=1/\la_{\min}$,
where $\la_{\min}$ is a minimal eigenvalue of $\fA$.
This establishes invertibility and the bound (5.16).

If in addition  $S$ is invertible and $\sigma$ is small 
 we may 
simply use the Neumann series to get invertibility of $\sigma \fA+S$. Namely, if $|\sigma|\le
\big(2\|\fA\| \|S^{-1}\|\big)^{-1}$ we
get
$(\sigma \fA+S)^{-1}=S^{-1}(I+\sum_{j=1}^\infty(-1)^j \sigma^j
(\fA S^{-1})^j)$ and the bound (5.17) is immediate.
\qed
\enddemo

\proclaim{Lemma 5.4} Let $\ell\ge 1$ be an odd integer, let $\Omega_1$ be the cone of real symmetric positive definite 
$\ell\times\ell$ matrices and let $\Omega_2$ be the set of all skew 
symmetric $\ell\times\ell$ matrices with rank $\ell-1$.

For $S\in\Omega_2$ choose a 
unit vector $e_S$ in the kernel of 
 $S$ and let $\pi_S$ be the orthogonal projection to the orthogonal complement of $e_S$.

Then for $A\in \Omega_1$, $S\in \Omega_2$, $\sigma\in \bbR$ we have 
$$
\det (\sigma A+S)=  \sigma
\inn{Ae_S}{e_S}
\det(\pi_S (\sigma A+S)\pi_S^*)
+ \sigma^2 F(A,S,\sigma)
\tag 5.18
$$
where 
$F$ is a smooth function 
on $\Omega_1\times \Omega_2\times \bbR$.
\endproclaim

\demo{Proof} Let $Q=Q(S)$ be an   orthogonal transformation with
$e_S^t Q=(0,\dots,1)$.
Then 
$$
Q^t(\sigma A+S)Q=\pmatrix \sigma A_0+S_0&\sigma a\\\sigma a^t& \sigma 
\eta
\endpmatrix
$$ where $S_0$ is a skew symmetric invertible $(\ell-1)\times(\ell-1)$
matrix, $A_0$ is positive definite, $a\in \bbR^{\ell-1}$ and $\eta=
\inn{Ae_S}{e_S}$. 
We apply Lemma 5.3 to $\sigma A_0 +S_0$ and  factor
$$\pmatrix \sigma A_0+S_0&\sigma a\\ \sigma a^t& \sigma\eta \endpmatrix=
\pmatrix
I&0\\ \sigma a^t(\sigma A_0+S_0)^{-1}&1\endpmatrix
\pmatrix
\sigma A_0+S_0&\sigma a\\0&\sigma\eta-\sigma^2 a^t(\sigma A_0+S_0)^{-1}a
\endpmatrix
$$
and conclude that
$$\det(\sigma A+S)=\det(\sigma A_0+S_0) 
\big(\sigma\eta-\sigma^2 a^t(\sigma A_0+S_0)^{-1}a\big).$$
The assertion follows since 
$\det(\sigma A_0+S_0)=\det(\pi_S (\sigma A+S)\pi_S^*)$.\qed
\enddemo

We now proceed to verify  the conditions (4.3-5) in \S4.
By Lemma 5.3 the determinant of
$\Phi_{(x,u),(y,v)}''$ can only vanish when 
$\sigma:=\sigma_{cr}\equiv x_d-e_d^t J_u y-u^t\La_d$ vanishes.
In this case the dimension of the kernel 
$\Phi_{(x,u),(y,v)}''$ is equal to the dimension of the kernel of $PJ_uP^t +E(B)$ with $B=\Gamma'(x'-y')$, thus equal to $1$. Thus $\rank(\Phi_{(x,u),(y,v)}'')\ge d+m-1$ everywhere.

In order to verify (4.4) let $V_L$ be a nonvanishing
 vector field  which is in the kernel of $Dp_L$ when
the mixed Hessian  (5.11)  becomes singular (i.e. when 
$x_d-e_d^t J_u y-u^t\La_d=0$).
Then
$$
V_L= \sum_{j=1}^{d-1} W_{L,j}\frac{\partial}{\partial y_j}
+
g_{L} \frac{\partial}{\partial y_d}+\sum_{i=1}^m h_{L,i} 
\frac{\partial}{\partial v_i},
\tag 5.19
$$
and with $A=\Gamma''(x'-y')$, 
we have 
$g_L=B^t W_L$ and
$$
(\sigma A+PJ_uP^t  +Be_d^t J_u P^t 
+PJ_u e_dB^t)W_L=0;
\tag 5.20
$$
moreover the functions $h_{L,i}$ are in the ideal generated by the 
$W_{L,j}$ (and the coefficients can be computed from (5.11)). 
To get a nontrivial kernel (when $\sigma=0$) we  must choose a nonvanishing vector
$W_L$ satisfying 
(5.20). Notice that then
$|e_d^t J_u P^t  W_L|$ is bounded below, by (5.14).
By Lemma 5.4 we have 
$$
V_L(\det \Phi_{(x,u),(y,v)}'')= (-1)^d 
F_1(x,y,u) e_d^t J_u P^t  W_L   +F_2(x,y,u,v) (x_d-e_d^tJ_u y-u^t\La_d)
$$
where $F_1$ and $F_2$ are smooth and  $F_1$ does not vanish.
Thus
$|V_L(\det \Phi_{(x,u),(y,v)}'')|\ge c$ on the zero set of
$\det \Phi_{(x,u),(y,v)}''$.

Next we consider the map $p_R$ and
let $V_R$ be a nonvanishing
 vector field  which is in the kernel of $Dp_R$ (or the cokernel of (5.11)) when
 $x_d-e_d^tJ_u y-u^t\La_d=0$.
Then
$$
V_R= \sum_{j=1}^{d-1} W_{R,j}\frac{\partial}{\partial x_j}
+
g_{R} \frac{\partial}{\partial x_d}+\sum_{i=1}^m h_{R,i} 
\frac{\partial}{\partial u_i}
$$
where by (5.11) the functions $h_{R,i}$ vanish when 
$x_d-e_d^tJ_u y-u^t\La_d=0$ and 
$$\align 
&W_R^t \big[\sigma A+PJ_u P^t +Be_d^tJ_u P^t ]+ g_R B^t=0
\\
&W_R^tPJ_u e_d-g_R=0;
\endalign
$$
thus since $A$ is symmetric and $J_u$ skew symmetric 
we have essentially the same equation for $W_L$ above, except that 
$J_u$ is replaced by $-J_u$:
$$
(\sigma A-PJ_u P^t -PJ_ue_dB^t -e_d^tJ_uP^t )W_R=0.
\tag 5.21
$$
Moreover  $g_R=e_d^t J_uP^t W_R$ does not vanish by (5.14).
As  $ x_d-e_d^tJ_u y-u^t\La_d$ does not depend on $x'$ we get
$$
V_R(\det \Phi_{(x,u),(y,v)}'')= 
\widetilde 
F_1(x,y,u) e_d^t J_u P^t  W_R   +\widetilde F_2(x,y,u,v) 
(x_d-e_d^tJ_u y-u^t\La_d)
$$
with smooth functions $\widetilde F_1$, $\widetilde F_2$ and nonvanishing
$\widetilde F_1$. Thus 
$|V_R(\det \Phi_{(x,u),(y,v)}'')|$ is bounded below on the
 zero set of
$\det \Phi_{(x,u),(y,v)}''$ and we have verified the statements analogous to (4.3-5).

\subheading{Proof of Claim 5.1, conclusion}
For small $l$ the bound is immediate from H\"ormander's 
standard $L^2$ estimate for nondegenerate oscillatory integrals
(\cite{8}, cf.  (5.12) and Lemma 5.3 above).
For large $l$ we can, by Lemma 5.4,  rewrite the amplitude  $\cE_j^l$ 
 as a finite sum 
$$\cE_j^l(x,y,u,v)= 2^{2jl} \sum_{|i|\le C} 
\zeta_1(2^{l+i}\det \Phi_{(x,u,y,v)}'') q_{l+i}(x,u,y,v)
$$ where the  $q_{l+i}$ are  compactly supported and smooth and satisfy
the estimates $\partial_{x,y,u,v}^\alpha q_{l+i}=O(2^{l\alpha})$.
 Since $2^l\le \la^{1/3}$  this type of blowup is covered by (4.2) and 
we can apply the estimate (4.6) and see that the operator
with kernel
$\la^{-j} \cE_j^l$ has $L^2$ operator norm
$\lc 2^{2jl}\la^{-j}\la^{-(d+m)/2} 2^{l/2}.$  This  implies our claim.

\subheading{Modifications for the proof of (3.16)}
By scaling we need to consider the operator of convolution with
$\partial_s K^{k,l}_s|_{s=1}$.

Let $\phi$ be as in (5.2)
and
$$
\align
&\rho(x',x_d,u,y,v,\sigma,\tau)=\frac{\partial}{\partial s}
\phi\big(\frac {x}{ s}, \frac {u}{s^2}, \frac{y}{s}, \frac{v}{s^2}, \sigma, \tau\big)\Big|_{s=1}
\\&=
\sigma \big(-x_d+y_d+(x'-y')\cdot\nabla_{x'}\Gamma(x'-y')\big)
+ 2\sum_{i=1}^m\tau_i(-u_i+v_i-x^tJ_iy)
+\sum_{i=1}^m\tau_i e_i^t \La(y-x).
\tag 5.22\endalign
$$

As before we set  $\lambda=2^k$  and observe that  our operator is 
a sum of an operator 
$\cG^{\la,l}$ with 
 Schwartz kernel  
$$
G^{\lambda,l}(x,u,y,v)= 
\lambda^{m+2}\iint  e^{i\la \phi(x,u,y,v,\sigma,\tau)}
\rho(x',x_d,u,y,v,\sigma,\tau) 
\chi_0(x,u,y,v) \eta_l( \sigma,\tau) d\sigma d\tau
$$
and an operator  which has 
similar properties as $H^{\lambda,l}$ above (thus satisfies 
estimates which are better than claimed in (3.16)).

We now need to  carry out the stationary phase calculations as before
for the kernel $\Cal F_{\la,1} \cG^{\lambda,l}$  (since the contribution from
$\Cal F_{\la,2}\cG^{\la,l}$  is again negligible). It has the form of (5.3), except that
$b_l$ is replaced by $\lambda c_l$ where $c_l$ is given by
$$
c_l(x,u,y,v,\theta)= 
b_l(x,u,y,v,z_d,w,\sigma,\tau)
\rho(x',z_d,w,y,v,\sigma,\tau)
. $$
Then by stationary phase the Schwartz kernel of $\cF_{\la,1}\cG^{\la,l}$
can be expanded as 
$$
\align
\la^{m+2}
\int e^{i\la \Psi(x,u,y,v,\theta)}c_{l}(x,u,y,v,\theta) d\theta&=
e^{i\la\Phi(x,u,y,v)} 
\sum_{j=0}^{N-1}\widetilde \cE_j^l(x,u,y,v)\lambda^{1-j}
+ \widetilde R_N^{\la,l}(x,u,y,v)
\tag 5.23
\endalign
$$
where again the error term $\widetilde R_N^{\la,l}$ is easy to handle for large $N$ and
$\widetilde \cE_j^\lambda$ 
is defined as in (5.8) but with $b_j$ replaced by
$c_j$.

In order to finish the proof of (3.16) it is now sufficient to establish that the operator 
$\cT^{\la,l}_j$ 
with
kernel $\la^{1-j}\widetilde \cE^l_j e^{i\la\Phi(x,u,y,v)}$ satisfies the bound
$$\|\cT^{\la,l}_j\|_{L^2\to L^2} \lc  \la^{1-(d+m)/2}2^{-l/2} 
(1+\|\La\|2^l).
\tag 5.24
$$ 
The differentiation in $s$ causes a blowup by not more than $\lambda$ and 
 by our previous analysis it follows that 
$$\|\cT^{\la,l}_j\|_{L^2\to L^2}\lc 2^{l/2}\la^{1-(d+m)/2} (2^{2l}\la^{-1})^j.
\tag 5.25
$$
If $j=1,2,\dots$ this estimate is sufficient for (5.24) since 
then $2^{l/2} (2^{2l}\la^{-1})^j\lc 2^{-l/2}$ 
by our restriction $2^l\le \la^{1/3}$.

This crude  estimate does not suffice for  the leading term  in the asymptotic expansion when $\|\La\|$ is small (or zero).

However note that when $\La=0$ 
 the coefficient of $\tau_i$ in (5.22)  
vanishes on the critical set where 
 $\theta=\theta_{\text{crit}}(x,u,y,v)$ since $\partial \Psi/\partial\tau=0$ on that set. We get
$$
\align
\rho(x',z_{d,\text{crit}}, w_{\text{crit}},y,v,\sigma_{\text{crit}},
\tau_{\text{crit}})=&
(x_d-e_d^tJ_u y-u^t \Lambda_d)
\big( (x'-y')\cdot\nabla_{x'}\Gamma(x'-y')-\Gamma(x'-y')\big)
\\&+2
\sum_{i=1}^m u_i 
\big(e_i^t\Lambda P^t (x'-y')+e_i^t\Lambda_d \Gamma(x'-y')\big).
\endalign
$$
Since $|x_d-e_d^tJ_u y-u^t\La_d|\approx 2^{-l}$ on the support of $c_l$ 
and since the coefficients of $u_i$ are $O(\|\La\|)$
 we now  gain an additional factor of $O(2^{-l}+\|\La\|)$ in the
 estimate  (5.25) for $j=0$  and thus establish (5.24) also for $j=0$.

\subheading{Modifications for  the proof of (3.17), (3.18)} The only reason for the 
modified definition  
(2.2.3) (replacing (2.2.2) for $l>k/3$) is the preservation of the symbol estimates (4.2), needed
for the validity of (4.6), (4.7). The estimation for 
$\widetilde K^k$ is exactly analogous to the estimation of $K^{k,l}$ when $l<k/3$, and the 
same statement 
applies to the $s$-derivatives. Only notational modifications are needed.

\head{\bf 6. Weak type (1,1) estimates}\endhead

We are now proving the weak type inequality (2.5). The proof of (2.6)
is omitted since it is exactly analogous.

We apply  standard Calder\'on-Zygmund arguments
 (with respect to nonisotropic families of balls 
on nilpotent Lie groups, see \cite{4}, \cite{17}). {\it Cf.}  also \cite{14}
and related papers on 
singular Radon transforms.

Let 
$$B_\delta=\{(x,u): |x|\le \delta, |u|\le \delta^2\}$$
and denote by $B_\delta^c$ its complement.

Since we have already checked  the $L^2$ bounds for the maximal function
it suffices to check the following H\"ormander type condition for $L^\infty(\bbR^+)$ valued kernels:
$$\sup_{\delta>0}\sup_{(y,v)\in B_\delta}\int_{B_{10\delta}^c}
\sup_{t>0}\big|
K^{k,l}_t\big((y,v)^{-1} (x,u)\big)-
K^{k,l}_t(x,u)\big| dx du \lc k 2^{k-l}(1+\|\La\|2^l)
$$
which follows from the two estimates
$$
\sup_{(y,v)\in B_\delta}\int_{B_{10\delta}^c}
\sup_{s\in[1,2]}\big|
K^{k,l}_{2^n s}\big((y,v)^{-1}(x,u)\big)-
K^{k,l}_{2^n s}(x,u)\big| dx du \lc \cases
& 2^{k-l}(1+\|\La\|2^l),
\\
&2^{k(m+2)}\min\{ 2^{-n}\delta, 2^n\delta^{-1}\}.\endcases
$$
Indeed we use the first bound for the $O(k)$ terms with $2^{-2k(m+1)}\le 2^{-n}\delta\le
2^{2k(m+1)}$ and the second bound for the remaining terms. We then sum the series in $n$.
Using scaling we see that the latter estimates are equivalent
to
$$\sup_{(y,v)\in B_r}\int_{B_{10r}^c}
\sup_{s\in[1,2]}\big|
K^{k,l}_{ s}\big(x-y,u-v+x^tJy)-
K^{k,l}_{ s}(x,u)\big| dx du \lc \cases
& 2^{k-l}(1+\|\La\|2^l),
\\
&2^{k(m+2)}\min\{ r^{-1}, r\}.\endcases
\tag 6.1
$$

Because of the support properties of the kernel 
 the integral on the left hand side is zero if $r\gg 1$.
Now assume that $r\lc 1$. 
 Since $|\nabla K^{k,l}_s(x,u)|\lc 2^{k(m+2)}$ the bound 
$2^{k(m+2)}r$  in (6.1) is immediate.
It remains to show that
$$\big\|\sup_{s\in[1,2]}|K^{k,l}_s|\big\|_1\lc 2^{k-l}(1+\|\La\|2^l),$$ 
and this follows from
$$\align
\big\|K^{k,l}\big\|_1&\lc 1,
\tag 6.2\\
\big\|\partial_s K^{k,l}_s\big\|_1&\lc 2^{k-l}(1+\|\La\|2^l).
\tag 6.3
\endalign
$$
By an integration by parts in $\sigma$, $\tau$ we see that
$$
|K^{k,l}(x,u)|\le C_N
\frac{2^{k-l}}{(1+2^{k-l}|
x_d-\Gamma(x')|)^N}
\frac{2^{km}}
{(1+2^k|u-\La x|)^N}
\tag 6.4
$$
from which (6.2) immediately follows.
Moreover from (5.22) one obtains
by the same argument 
$|\partial_s K^{k,l}_s(x,u)|$ 
is bounded by  $ C_N' 2^{k-l}(1+\|\La\|2^l)$ times the right hand side of 
(6.4). Consequently we obtain (6.3). This finishes the proof of the weak type inequality
(2.5).
\qed


\head{ \bf 7. Appendix}\endhead

In this section we give the example  of a two-step nilpotent Lie group 
$G$,
with $10$-dimensional Lie algebra, which satisfies 
 the nondegeneracy condition but which is  not isomorphic to a group  of 
Heisenberg type.

For $\mu=(\mu_1,\mu_2)\in \bbR^2$ let
$$E_\mu=\pmatrix 
\mu_1&0&0&-\mu_2
\\
\mu_2&\mu_1&0&0
\\
0&\mu_2&\mu_1&0
\\
0&0&\mu_2&\mu_1
\endpmatrix$$ and define the $8\times 8$ matrix
$$J_\mu=\pmatrix 0&E_\mu\\-E_\mu^t&0\endpmatrix;$$
then $$\det J_\mu=(\mu_1^4+\mu_2^4)^2.\tag 7.1
$$
Let $\fg$ be the Lie algebra which is $\bbR^8\oplus \bbR^2$ as a vector space, with Lie bracket
$$
[X+U,Y+V]=0 + (X^tJ_{(1,0)}Y,  X^tJ_{(0,1)}Y).
$$
By (7.1) the group identified with $\fg$ satisfies our nondegeneracy condition.
 We now prove by contradiction 
 that $\fg$ is not isomorphic to a Heisenberg-type Lie algebra.

Assume that there is a Lie algebra isomorphism
$\alpha: \tfg\to \fg$ where $\tfg$ is a Heisenberg-type algebra.
Then $\tfg=\fw\oplus \fz$ where $\fz$ is the center and $\alpha$ is 
a linear isomorphism from $\fz$ to $\Bbb R^2$.

Now with respect to orthonormal bases $u_1,\dots, u_8$ on $\fw$ 
and $u_9, u_{10}$ on $\fz$ and $ e_1,\dots, e_8$ on $\bbR^8$ and 
$e_9,
 e_{10}$ on $\bbR^2$ the map $\alpha $ is given by the $10\times 10$
 matrix
$$\pmatrix A&0\\L&B\endpmatrix$$ where $A$ is an invertible $8\times 8$ matrix and $B$ an invertible  $2\times 2$ matrix. 

Now let $X=\sum_{i=1}^8 x_i u_i$, $Y=\sum_{i=1}^8 y_i u_i$, 
and express $\omega\in \fz^*$ in terms of the dual basis as 
$\omega=w_1u_9^*+w_2 u_{10}^*$.
Then, since $\tfg$ is of Heisenberg type we have 
$\omega([X,Y])= x^t \widetilde J_w y$
 with $\widetilde J_w^2=-(w_1^2+w_2^2) I$; in particular
$$|\det \widetilde J_w|=(w_1^2+w_2^2)^4.\tag 7.2$$
Now if $\omega=\alpha^t \mu$  
(thus $B^t\mu=(w_1,w_2)^t$) then
$$x^t \widetilde J_{B^t\mu} y=
\omega([X,Y])= 
(\alpha^t)^{-1} \omega(\alpha[X,Y])= 
\inn{\mu}{[\alpha X,\alpha Y]}
= (Ax)^tJ_\mu (Ay)
$$
so that $A^t J_\mu A=\widetilde J_{B^t\mu}$
and therefore 
$$
\det \widetilde J_{B^t\mu}=(\det A)^2 \det J_\mu.
$$
Thus  by (7.1) and (7.2) we obtain $|B^t\mu|^8= (\det A)^2 (\mu_1^4+\mu_2^4)^2$ and 
therefore, if  $(a,b)$ and $(c,d)$ are 
the rows of the matrix $|\det A|^{-1/4}B^t$,
$$
\mu_1^4+\mu_2^4= 
\big( (a\mu_1+b\mu_2)^2+(c\mu_1+d\mu_2)^2\big)^2,
$$
for all $\mu\in \bbR^2$. 
Thus
$$\mu_1^4+\mu_2^4= 
\big( (a^2+c^2)\mu_1^2 +(b^2+d^2)\mu_2^2+ 2(ab+cd) \mu_1\mu_2)\big)^2
$$
for all $\mu\in \bbR^2$.  This implies $a^2+c^2=b^2+d^2=1$ and setting
$\rho=ab+cd$  we obtain after a little algebra that
$$(4\rho^2 +2)\mu_1\mu_2+ 4\rho(\mu_1^2+\mu_2^2)=0$$
for all  $\mu\in \bbR^2$. This implies both $2\rho^2+1=0$ and $\rho=0$,
thus a contradiction.\qed

\enddocument